\documentclass{article}

\usepackage{amsmath,amsfonts,amssymb,amsthm}
\font \maius =cmcsc10
\newtheorem{lemma}{Lemma}[section]
\newtheorem{theorem}[lemma]{Theorem}
\newtheorem{prop}[lemma]{Proposition}
\newtheorem{corollary}[lemma]{Corollary}
\theoremstyle{definition}
\newtheorem*{dfn}{Definition}

\newtheorem{remark}[lemma]{Remark}

\theoremstyle{plain}

\newcommand{\Q}{\mathbb{Q}}\newcommand{\R}{\mathbb{R}}
\newcommand{\C}{\mathbf{C}}

\newcommand{\LL}{\mathbb{L}}
\newcommand{\Z}{\mathbb{Z}}
\newcommand{\F}{\mathbf{F}}
\newcommand{\HH}{\mathbf{H}}

\newcommand{\B}{\mathbf{B}}
\newcommand{\N}{\mathbb{N}}
\newcommand{\m}{\ensuremath{\mathfrak m}}

\newcommand{\x}{\ensuremath{\textbf{x}}}
\newcommand{\y}{\ensuremath{\textbf{y}}}
\newcommand{\z}{\ensuremath{\textbf{z}}}
\newcommand{\e}{\ensuremath{\textbf{e}}}
\newcommand{\p}{\ensuremath{\textbf{p}}}

 \DeclareMathOperator{\Der}{Der}

\DeclareMathOperator{\Ad}{Ad}
\DeclareMathOperator{\ad}{ad}
 \DeclareMathOperator{\End}{End}
\DeclareMathOperator{\GL}{GL}

\DeclareMathOperator{\Aut}{Aut} \DeclareMathOperator{\car}{char}

\title{On linearity of finitely generated $R$-analytic groups
\thanks{This work has been partially  supported by the FEDER, the
MCYT Grants BFM2001-0201, BFM2001-0180  and the Ram\'on y Cajal
Program.} }
\author{ A. Jaikin-Zapirain\\
\footnotesize Departamento de Matem\'aticas
\footnotesize Facultad de Ciencias\\
\footnotesize Universidad Aut\'onoma de Madrid}
\date{}

\begin{document}

\maketitle

\begin{abstract}
We prove that if $R$ is a commutative Noetherian local pro-$p$
domain of characteristic 0 then every finitely generated
$R$-standard group is $R$-linear.
\end{abstract}

\section{Introduction}\label{intr}
Let $R$ be a commutative Noetherian local pro-$p$ domain and $\m$
its maximal ideal.
 The concept of an $R$-analytic group is
defined in \cite[Chapter 13]{DDMS}, where it is shown that if $R$
satisfies some additional technical conditions, then every such
group contains an open subgroup which is \textbf{$R$-standard}. To
recall what this means, let $G$ be an $R$-standard group. Then the
underlying set of $G$ may be ``identified'' with the cartesian
product $(\m^l)^{(d)}$ of $d$ copies of $\m^l$, for some $l \in
\N$. The number
 $d\ge 0$ is the \textbf{dimension} of $G$ and $l>0$ is the \textbf{level} of $G$.
 The group operation is given by a \textbf{formal group law},
i.e.\ a $d$-tuple $\mathbf{F}=(F_1,\ldots,F_d)$ of power series over
$R$ in $2d$ variables, as follows: for all $\x,\y \in G =
(\m^l)^{(d)}$ we have
$$
\x \cdot \y = (F_1(\x,\y),\ldots,F_d(\x,\y)).
$$
The neutral element of $G$ is $\e=(0,\ldots,0)$. Without loss of
generality   we will always  assume that the level of $G$ is 1. We
shall write $G(I)$ for $(I)^{(d)}\subset G$ and $G_i$ will denote
$G(\m^i)$. Hence $G=G_1=G(\m)$.

The $\Z_p$-analytic pro-$p$ groups are well-understood (see
\cite{La,DDMS})  and they are linear over $\Z_p$. It was
conjectured
 that $R$-analytic pro-$p$ groups are $R$-linear. Since compact $\Z_p$-analytic groups are finitely generated, it is natural to
 consider first finitely generated $R$-analytic pro-$p$ groups.
  In \cite{Ja} it was proved that  just infinite $R$-analytic groups are $R$-linear.
   As an $R$-analytic group contains an open $R$-standard subgroup,
   we can restrict our attention to $R$-standard groups. In \cite{CadS} the linearity of  $\Z_p[[t ] ]$-perfect groups is
    shown. Recall that an $R$-standard group of level $l$ is called $R$-perfect if $[G,G
    ]=G_{2l}$. Note that $R$-perfect  groups are finitely
    generated.
     In this work  we prove:
\begin{theorem} \label{main}
Let $R$ be commutative Noetherian local pro-$p$ domain of
characteristic 0  and $G$ a finitely generated $R$-standard group.
Then $G$ is $R$-linear.
\end{theorem}
Note that if a pro-$p$ group is finitely generated and
$R$-linear, then it is a closed subgroup of $\GL_n(R)$ for some
$n$. In the following we  will use the notion $R$-linear only for
closed subgroups of $\GL_n(R)$.

 First at all in Section \ref{t-linear} we prove some new results
 about t-linear pro-$p$-groups.
   In Section \ref{liealgebra} we shall describe
 the Lie algebra $\LL=\LL(G)$ of an $R$-standard group $G$. The typical definition is based on
  the formal group law of $G$. We shall introduce another definition coming from the theory of
  algebraic groups and define $\LL$ as the set of left invariant derivations of $A=R[[x_1,\cdots x_d ] ]$.
  We shall prove that these two definitions are equivalent (it is a folklore result,  however,
  I do not know any reference for it in the literature).
In Section \ref{useBCHF}, assuming that the characteristic of $R$
is 0, we use the BCHF to define on $\p\LL$ (here $\p=\p(p)=4$ if
$p=2$ and $\p=p$ if $p$ is an odd prime) a group structure. We
show that the obtained group is isomorphic to $G(\p R)$.   In
Section \ref{nilprad} we show that if an $R$-standard group is
finitely generated, then the radical of $\LL(G)$ is nilpotent. It
will permit us use the Weigel result about linearity of some Lie
rings which we describe in Section \ref{prel}. In Section
\ref{final} we finish the proof of Theorem \ref{main}.

\section{Some results about t-linear pro-$p$ groups}\label{t-linear}
Recall the definition of t-linear pro-$p$ groups from \cite{Ja}.
\begin{dfn}
Let $G$  be a pro-$p$ group. We shall say that $G$ is {\em
t-linear} if it is a closed subgroup of $\GL_n(A)$ for some
commutative profinite ring $A$.
\end{dfn} In \cite[Theorem 4.1]{Ja} it was shown that if $G$ is a finitely generated t-linear pro-$p$ group, then
 $G$ is linear over some
commutative Noetherian local pro-$p$ ring. In this section we
extend this result. If $R$ is a ring, we denote by $K\dim R$ the
Krull
 dimension of $R$.
\begin{theorem}\label{tlinear} Let $G$ be a finitely generated  pro-$p$
group and suppose that $G$ is linear over some commutative
Noetherian local pro-$p$ domain $R$. Then we have

1. If $\car R=0$ and $K\dim R> 2$, then $G$ is linear over every
commutative Noetherian local pro-$p$ domain of characteristic zero
and Krull dimension greater than 2.

2. If $\car R=0$ and $K\dim R\le 2$, then $G$ is linear over every
commutative Noetherian local pro-$p$  domain of characteristic
zero and Krull dimension greater or equal than of $R$.

3.  If $\car R=p$ and  $K\dim R\ge 2$, then $G$ is linear over
every commutative Noetherian local pro-$p$ domain $R$ of
characteristic $p$ and
 Krull dimension greater than 1.
 \end{theorem}
 We see that the last theorem reduces the study of
linear over pro-$p$ domains pro-$p$ groups to study of $R$-linear
pro-$p$ groups,
 where $R=\Z_p[[t_1,t_2]]$ or $R=\F_p[[t_1,t_2]]$.

 One of the main step in the proof of the previous theorem is the
 following proposition. We also will use it in the proof of our
 main result.
 \begin{theorem}\label{torsionfreemodule}
 Let $R$ be a commutative
Noetherian local pro-$p$ domain and $W$ a finitely generated
$R$-torsion-free $R$-module. Suppose $G\le\Aut_R(W)$. Then there
exists a commutative Noetherian local pro-$p$ domain $S$,
satisfying $\car S=\car R$ and $K\dim S=K\dim R$, such that $G$ is
$S$-linear. Moreover, if $R$ is regular, then $S=R$.
\end{theorem}
First we need some auxiliary results. In the following $R$ is
always a commutative Noetherian local pro-$p$ domain.
\begin{lemma}\label{kr}
Let $a\ne 0,r\in R$ and $T=R[[t ] ]/(a t-r)$. Then the Krull
dimensions of $R$ and $T$ are the same.
\end{lemma}
\begin{proof} First note that $K\dim R[[t ] ]=K\dim R+1$ and
since $R[[t ] ]$ is domain $K\dim R[[t ] ]/(a t-r)$ is strictly
less than $K\dim R[[t ] ]$. Hence $K\dim R[[t ] ]/(a t-r)\le K\dim
R$.

On the other hand, $K\dim T$ is equal to the number of elements in
a system of parameters of $T$ (see \cite[p.27]{Na}), and this
number is at least $K\dim R[[t ] ]-1$. This implies $K\dim R[[t ]
]/(a t-r)\ge K\dim R$.
\end{proof}

Recall that a {\bf (multiplicative non-archimedean) valuation} of
a field $D$ is a mapping $v\colon D\to \R_{\ge 0}$ such that for
$a,b\in D$.
\begin{enumerate}
\item
$v(a)=0$ if and only if $a=0$.
\item
$v(ab)=v(a)v(b)$.
\item
$v(a+b)\le \max (v(a),v(b))$.
\end{enumerate}
We need the following proposition:
\begin{prop}(\cite[Proposition 11.9]{Na})\label{valuation} Let $R$ be a subring of
a field $D$ and $\m$ a non-trivial ideal of $R$. Then there exists
a valuation $v$ of $D$ such that $v(r)\le 1$ for every $r\in R$
and $v(m)<1$ for every $m\in \m$.
\end{prop}
\begin{lemma}\label{intclous}
Let $S$ be the integrally closure of $R$. Then $S$ is also a
commutative Noetherian local pro-$p$ ring of same characteristic
and  same Krull dimension as $R$.
\end{lemma}
\begin{proof} First we want to see that $S$ is local. Let $D$ be the  quotient field of $R$ and $v$ a valuation from
Proposition \ref{valuation}. If $s\in S$ then $v( s)\le 1$,
because $s$ is integral over $R$. In order to see that $S$ is
local it is enough to show that if $v(s)=1$ then $s$ is invertible
in $S$. Let $f(t)=\sum_{i=0}^{n}a_it^i$ be a monic irreducible
over $R$ polynomial such that $f(s)=0$. Since $v(s)=1$, there
exists $a_i$, $i<n$, such that $a_i\notin \m$.  Therefore, since
$R$ is a Henselian ring (see \cite[Theorem 30.3]{Na}), we have
that $a_0\notin \m$. Hence $s^{-1}\in S$.

Finally, by \cite[Theorem 32.1]{Na}, $S$ is a finite extension of
$R$, whence their Krull dimensions coincide.
\end{proof}
\begin{theorem}\label{embed} Let $R$ be a commutative
Noetherian local pro-$p$ domain and $a\in R$. Then there  are a
commutative Noetherian local pro-$p$ domain $S$ and an injective
homomorphism $\phi\colon R\to S$ such that
\begin{enumerate}
\item
$\phi(\m^k)\subseteq \phi(a )S$ for some $k\in \N$;
\item
$S$  is integrally closed and its Krull dimension is the same as
of $R$. \end{enumerate}
Moreover, if $R$ is regular, then $S=R$.
\end{theorem}
\begin{proof}
Let $D$ be the quotient field of $R$ and $v$ a valuation from
Proposition \ref{valuation}.  Let $\Omega$ be a completion of $D$
respect to $v$. Since $R$ is Noetherian, there exists $s=\max \{
v(r)|r\in \m\}<1$ Let $k$ be such that $s^k/v(a)<1$ and put
$S_1=R[[\m^k/a ] ]\subset \Omega$. It is clear that $S_1$ is a
local pro-$p$ ring. Moreover, if $t_1,\ldots,t_m$ are generators
of $\m^k$ as $R$-module, then $S_1=R[[t_1/a,\ldots,t_m/a ] ]$.
Hence $S_1$ is Noetherian. Applying several times Lemma \ref{kr},
we obtain that its Krull dimension is the same as of $R$. Finally,
let $S$ be the integral closure of $S_1$. Theorem follows from
Lemma \ref{intclous}.

If $R=A[[t_1,\ldots,t_l] ]$, where $A$ is equal to $\F_q$ or to a
finite extension of $\Z_p$, then the homomorphism $\phi\colon R\to
R$ is defined by means of $\phi(t_i)=a t_i$.
\end{proof}
\begin{proof}[Proof of Theorem \ref{torsionfreemodule}]
Let $D$ be the field of quotients of $R$. Consider the $D$-module
$M=D\otimes_R W$. Since $W$ is $R$-torsion-free, we can see $W$ as
an $R$-submodule of $M$.

Let $m_1,\cdots, m_r$ be a $D$-basis of $M$ lying in $W$. Put
$N=\sum Rm_i$. It is clear that $N$
 is a free $R$-module. Let $a\ne 0$ be such that $aW\le N$. By the
 previous result, there are a commutative Noetherian local pro-$p$ ring $S$ and an
injective homomorphism $\phi\colon R\to S$ such that
$\phi(\m^k)\subseteq \phi(a) S$ for some $k\in \N$. Moreover, if
$R$ is regular, then $S=R$.

Put $W_1=\m^k W$ and $G_1=\{g\in G|gw\equiv w \pmod {W_1}
\textrm{\ for every \ }w\in W\}$. It is clear that $G_1$  is of
finite index in $G$ and $G$ acts faithfully on $W_1$. Define
$L=S\otimes _R N$. We have $L$ is a free $S$-module.

Now we  embed $W_1$ in $L$ in the following way. Let $w\in W_1$.
Then $w=\sum_{i=1}^r a^{-1}k_im_i$, where $k_i\in \m^k$. Define
$\psi(w)=\sum_{i=1}^r \phi(a)^{-1}\phi(k_i)\otimes m_i$  (since
$S$ is domain and $\phi(\m^k)\subseteq \phi(a) S$, we can speak
about $\phi(a)^{-1}\phi(k_i)\in S$). The map $\psi$ is an
$R$-homomorphism. So we can see $W_1$ as $R$-submodule of $L$.

Now, we explain how we can extend the action of $G_1$ on $L$. Let
$g\in G_1$ and $l=\sum_{i=1}^r s_im_i$. Define $gl= \sum_{i=1}^r
s_igm_i$. Note that, since $g\in G_1$, $gm_i\in W_1+N$, so the
definition is correct. This action gives   an embedding $G_1\le
\Aut_S(L)\cong \GL_r(S)$. Since
 $G_1$ is of finite index in $G$, $G$ is also $S$-linear.
 \end{proof}
\begin{proof}[Proof of Theorem \ref{tlinear}]By the structure theorem of complete local rings (see
\cite[Corollary 31.6]{Na}),  $R$ is  a finite extension of a
regular ring $T=\Z_p[[t_1,\cdots,t_k]]$ or
$T=\F_p[[t_1,\cdots,t_k]]$ for some $k$. Hence, $G\in \Aut_T(R^n)$
for some $n$, and Theorem \ref{torsionfreemodule} implies that $G$
is $T$-linear. Note that the Krull dimensions of $R$ and $T$ are
the same.

1.  If $\car T=0$ and $K\dim T> 2$. Then by Remark VII.10.4 of
\cite{ZaSa}, we can embed $T=\Z_p[[t_1,\cdots ,t_k ]]$ into
$\Z_p[[s_1,s_2]]$. On the other hand $\Z_p[[s_1,s_2]]$ can be
embeded into every commutative Noetherian local pro-$p$ domain $S$
of characteristic zero and Krull dimension greater than 2, and so
$G$ is also $S$-linear.

The proofs of 2. and 3. follow  the same ideas.
\end{proof}

\section{Lie algebra of an  $R$-standard group}\label{liealgebra}
We use the notation of Section \ref{intr}. So $G$ is an
$R$-standard group of level 1. The law $\F$ can be written in the
form
$$\F(\x,\y)=\x+\y+\B(\x,\y)+O^\prime(3),$$
 where $\B$ is the sum of the all polynomials in $\x$ and
 $\y$ of degree 2 and the expression $O^\prime(n)$ stands for any power series in which every term
 has total degree at least $n$ and has degree at least 1 in each variable. We know, see, for example,
 \cite[p.26]{newhor}, that if $\C(\x,\y)=\B(\x,\y)-\B(\y,\x)$, then $(R^{(d)},+,\mathbf{C})$ is a Lie $R$-algebra,
 and we shall denote this Lie algebra by $\LL=\LL(G)$.

Now, let $A=R[[x_1,\ldots,x_d] ]$. Since $G$ is identified with
$\m^{(d)}$, $A$ can be considered as a subring of the ring of
functions from $G$ to $R$. Note that since $R$ is domain, two
different elements from $A$ give us two different functions.
Define two actions of $G$ on $A$ via left and right translation:
$$(\lambda_{\x}f)(\y)=f(\x^{-1}\y),\ (\rho_{\x}f)(\y)=f(\y\x),\textrm{\ where \ } f\in A,\x\in G,\y \in G.$$
Since the multiplication in $G$ is given by an  analytic function,
we have that if $f\in A$ then $\lambda_{\x}f$ and $\rho_{\x}f$
also belong to $A$.

  The  bracket of two  $R$-derivations of $A$ is again a derivation. Therefore, $\Der A$ is  a Lie algebra.
  So is the subspace of left invariant derivations
  $\Der_lA=\{w\in \Der A\mid w\lambda_{\x}=\lambda_{\x}w \textrm{\ for all\ }\x\in G\} $,
  since the bracket of two derivations which commute with $\lambda_{\x}$ obviously does likewise.
The next theorem is the main result of this section:
 \begin{theorem}\label{equiv}
The Lie $R$-algebras $\LL$ and $\Der_lA$ are isomorphic.
\end{theorem}
Before the proof of the theorem we need to do some preliminary work.
\begin{lemma}
Let $w_1,w_2\in \Der_lA$ and $f\in A$. Then we have
\begin{eqnarray*}
w_1(f)(\y) & =& \sum_{i=1}^d\left. \frac{\partial f(\y\x)}{\partial x_i}\right |_{\x=\e}w_1(x_i)(\e),\\
w_2(w_1(f))(\z)& = & \sum_{i,j=1}^d\left. \frac{\partial^2f(\z\y\x)}{\partial x_i\partial y_j}\right |_{(\x,\y)=(\e,\e)}w_1(x_i)(\e)w_2(x_j)(\e).
\end{eqnarray*}
\end{lemma}
\begin{proof}
We prove only the first equality, because the second one is obtained applying two times the  first.
\begin{eqnarray*}
w_1(f)(\y)& = &(\lambda_{\y^{-1}} w_1(f(\x)))(\e)= w_1(\lambda_{\y^{-1}}f(\x))(\e)=w_1(f(\y\x))(\e)\\
& = & \sum_{i=1}^d\left. \frac{\partial f(\y\x)}{\partial x_i}\right |_{\x=\e}w_1(x_i)(\e).
\end{eqnarray*}
 \end{proof}
\begin{proof}[Proof of Theorem \ref{equiv}]
Define an $R$-homomorphism $\phi\colon \Der_l(A)\to \LL$ as follows
$$\phi(w)=(w(x_1)(\e),\ldots,w(x_d)(\e)), w\in \Der_l(A).$$
First we will  show that $\phi$ is a bijective map. Fix $(a_1,\ldots,a_d)\in \LL$ and define  $w\in \End_{R}(A)$ by means of
$$w(f)(\y)=\sum_{i=1}^d\left. \frac{\partial f(\y\x)}{\partial x_i}\right |_{\x=\e}a_i, f\in A.$$
If $f_1,f_2\in A$, then
\begin{eqnarray*}
w(f_1f_2)(\y) & = & \sum_{i=1}^d\left. \frac{\partial f_1(\y\x)f_2(\y\x)}{\partial x_i}\right |_{\x=\e}a_i\\
&=& \sum_{i=1}^d\left. \frac{\partial f_1(\y\x)}{\partial x_i}\right |_{\x=\e}f_2(\y)a_i+ \sum_{i=1}^df_1(\y) \left. \frac{\partial f_2(\y\x)}{\partial x_i}\right |_{\x=\e}a_i\\
& = & w(f_1)(\y)f_2(\y)+f_1(\y)w(f_2)(\y).
\end{eqnarray*}
This implies that $w$ is a derivation of $A$. Now, if $\z\in G$, then
\begin{eqnarray*}
w(\lambda_{\z}f)(\y)&=& \sum_{i=1}^d\left. \frac{\partial (\lambda_{\z}f)(\y\x)}{\partial x_i}\right |_{\x=\e}a_i\\
&=&w(f)(\z^{-1}\y)=(\lambda_{\z}w(f))(\y).
\end{eqnarray*}
Hence, we obtain that $w$ is really a left  invariant derivation.
Define the constructed map from $\LL$ to $\Der_l(A)$ by $\psi$.
Note that if $w\in \Der_l(A)$, then  we have
\begin{eqnarray*}
\psi(\phi(w))(f)(\y)&=&\psi(w(x_1)(\e),\ldots,w(x_d)(\e))(f)(\y)\\
&=& \sum_{i=1}^d\left. \frac{\partial f(\y\x)}{\partial x_i}\right |_{\x=\e}w(x_i)(\e)\\
&=&w(\lambda_{\y^{-1}} f)(\e)=\lambda_{\y^{-1}} w(f)(\e)=w(f)(\y).
\end{eqnarray*}
On the other hand if $a=(a_1,\ldots,a_d)\in \LL$, then
\begin{equation*}
\phi(\psi(a))= (b_1,\ldots,b_d),
\end{equation*}
where $b_k=\sum_{i=1}^d\left. \frac{\partial F_k(\y,\x)}{\partial x_i}\right |_{\x=\e}a_i=a_k$. Hence $\psi(\phi(a)=a$. We conclude that $\phi$ is a bijection.

 We shall see now that $\phi$ is also a homomorphism of Lie rings.

Let $w_1,w_2\in \Der_l(A)$. From the previous lemma, we obtain that
$$\phi([w_1,w_2])=(c_1,\ldots,c_d),$$
where $$c_k= \sum_{i,j=1}^d\left. \frac{\partial^2F_k(\y,\x)}{\partial x_i\partial y_j}\right |_{(\x,\y)=(\e,\e)}(w_2(x_i)(\e)w_1(x_j)(\e)-w_1(x_i)(\e)w_2(x_j)(\e)).$$
We conclude that $\phi([w_1,w_2])=\C(\phi(w_1),\phi(w_2)).$
\end{proof}
\begin{remark}\label{ident}
Note that the previous proof  also gives an identification of
$\Der_l(A)$ with $(I/I^2)^*=\End_R(I/I^2,R)$. The derivation $w$
is identified with the map $f\to w(f)(\e)$. \end{remark} In the
rest of the work we shall use the letter $\LL$ for $\Der_l(A)$.

\section{An application of the Baker-Campbell-Hausdorff Formula}\label{useBCHF}
The Baker-Campbell-Hausdorff Formula (BCHF) is $H(x_1,x_2)=\log(e^{x_1}e^{x_2})$ regarded as a formal power series in two non-commuting variables. Equivalently, this is a formal power series $H(x_1,x_2)$  such that
$$e^{H(x_1,x_2)}=e^{x_1}e^{x_2}.$$
The homogeneous component of $H(x_1,x_2)$ of degree $n$  is
denoted by $H_n(x_1,x_2)$, so that $H(x_1,x_2)=\sum_{n=1}^\infty
H_n(x_1,x_2)$. The main fact about the BCHF is that $H_n(x_1,x_2)$
is a Lie word in $x_1$ and $x_2$ (see, for example, \cite[Theorem
9.11]{Kh}). Let $\{e_k\}$ be  a $\Z$-basis of the free Lie algebra
generated by $x_1$ and $x_2$, consisting of simple commutators.
Then we can express $H_n$ as $H_n=\sum \lambda_{k,n}e_k$, for some
$\lambda_{k,n}\in \Q$. We need the following fact about the
coefficients $\lambda_{k,n}$ (\cite[Proposition II.8.1]{Bou}):
\begin{equation}\label{lamb}
v_p(\lambda_{k,n})\ge -(n-1)/(p-1),
\end{equation}
where $v_p(ap^s)=s$ if $(a,p)=1$.
\begin{theorem}\label{liegroup}
Let $M$ be a Lie $\Z_p$-algebra without $\Z_p$-torsion and suppose that $\cap_i p^iM=0$. If $M$ is complete in the topology induced by the filtration  $p^iM$, then
\begin{enumerate}
\item
If $a,b\in \p M$, then $H_n(a,b)\in \p M$ and $\lim_{n\to \infty}H_n(a,b)=0$ (in particular, we can compute $H(a,b)$);
\item
$(\p M, H)$ is a group.
\end{enumerate}
\end{theorem}
\begin{proof}
The first proposition of the theorem follows directly from
the formula (\ref{lamb}).

Since $H(a,0)=H(0,a)=a$ and $H(a,-a)=0$ for any $a\in \p M$, then in order to prove the second statement, we only need to show that the operation $H(x_1,x_2)$ is associative.

Let $F$ be the $\Q_p$-algebra of formal power series in the non-commuting variables $x_1,x_2,x_3$. Then we have the following equalities in $ F$:
\begin{eqnarray*}
H(H(x_1,x_2),x_3)& =& \log (e^{H(x_1,x_2)}e^{x_3})=\log(e^{x_1}e^{x_2}e^{x_3})\\
& =& \log(e^{x_1}e^{H(x_2,x_3)})=H(x_1,H(x_2,x_3)).
\end{eqnarray*}
In particular, we obtain
\begin{equation}\label{ass}
H(H(\p x_1,\p x_2),\p x_3)=H(\p x_1,H(\p x_2,\p x_ 3)).
\end{equation}
Let $L$ be  the  Lie  $\Z_p$-subalgebra of $F^{(-)}$, generated by
$x_1,x_2,x_3$. It is clear that $L$ is a free Lie $\Z_p$-algebra.
Let$\{e_k, k\in \N\}$ be a $\Z_p$-basis of $L$ and $ \bar L$ be
the Lie $\Z_p$-subalgebra of $F^{(-)}$, consisting from the formal
power series $\sum \lambda_ke_k$ with $\lambda_k\in \Z_p$ and
$\lim_{k\to \infty} (\lambda_k)=+\infty$. We have that $\bar L$ is
 the completion of $L$ in the topology induced by the
filtration $p^iL$.  By (\ref{lamb}), if $y_1,y_2\in \p\bar L$,
then $H(\p y_1,\p y_2)\in \p \bar L$.

Now, let $a,b,c \in M$ and $\phi\colon L\to M$ be a Lie $\Z_p$-algebra homomorphism defined by means of
$\phi(x_1)=a,\phi(x_2)=b,\phi(x_3)=c$. Since $M$ is complete in the topology induced by the filtration  $p^iM$,
this homomorphism can be extended to $\phi\colon \bar L\to M$. Furthermore, $\phi$ is a continuous map. Then
\begin{eqnarray*}
H(H(\p a,\p b),\p c)& =& H(H(\phi(\p x_1),\phi(\p x_2)),\phi(\p x_3))\\
& = &\phi(H(H(\p x_1,\p x_2),\p x_3)) = \phi(H(\p x_1,H(\p x_2,\p x_ 3)))\\
& = & H(\p a,H(\p b,\p c)).
\end{eqnarray*}
We conclude that the operation $H(x_1,x_2)$ is associative and, so, $(\p M,H)$ is a group.
\end{proof}
We will denote the group $(\p M,H)$ by $\Gamma(\p M)$.
\begin{theorem}\label{explog}
Let $D$ be an associative $\Z_p$-algebra without $\Z_p$-torsion and suppose that $\cap_i p^iD=0$. Assume that $D$ is complete in the topology induced by the filtration  $p^iD$. If $M$ is a closed Lie $\Z_p$-subalgebra of $D^{(-)}$,then
\begin{enumerate}
\item
If $a\in \p D$, then $a^n/n!\in \p D$ for $n\ge 1$ and $\lim_{n\to \infty}a^n/n!=0$ (in particular, we can compute $e^a\in 1+\p D$ and $\log(1+a)\in \p D$);
\item
If $a\in \p D$, then  $\log(e^a)=a$ and $e^{\log(1+a)}=1+a$ (in particular, $e^a=1$ if and only if $a=0$);
\item
$\{e^a|a\in \p M\}$ is a group (with multiplication of $D$) isomorphic to $\Gamma(\p M)$.
\end{enumerate}
\end{theorem}
\begin{proof}
The first statement follows from the fact that $v_p(n!)\le (n-1)/(p-1)$ (see \cite[Lemma II.8.1]{Bou}).

Let  $P$ be the subring of $\Z_p[[t] ]$, consisting of the series
$$\sum_i\alpha_it^i \textrm{\ with\ } \lim_{i\to  \infty}v_p(\alpha_i)=+\infty.$$
Suppose  $f\in t\p P$. Since $v_p(n!)\le (n-1)/(p-1)$,  $e^{f}\in P$ and $\log(1+f)\in P$. Let $a=\p b$. Define a homomorphism of $\Z_p$-algebras $\phi:\Z_p[t ]\to D$, by means of $\phi(t)=b$. Since $P$ is isomorphic to the completion of $\Z_p[t ]$ in the topology induced by $p^k\Z_p[t ]$, we can extend $\phi$ on $P$. Note that in $P$, the equality $\p t=\log(e^{\p t})$ holds. Since $\phi$ is  continuous, we have
$$a=\phi(\p t)=\phi(\log(e^{\p t}))=\log(e^{a}).$$
 Analogically, $e^{\log(1+a)}=a$. This proves the second proposition.

Since $M$ satisfies the hypothesis of the previous theorem, in order to prove the third statement we have to show that $e^ae^b=e^{H(a,b)}$, for any $a,b\in \p M$. The proof of this equality is analogical of the proof of Theorem \ref{liegroup}, and we omit it.
\end{proof}

\begin{lemma}\label{derivation}
Let $D$ be an associative $\Z_p$-algebra  without $\Z_p$-torsion
and suppose that $\cap_i p^iD=0$. Assume that $D$ is complete in
the topology induced by the filtration  $p^iD$. Let $\phi$ be a
$\Z_p$-automorphism of $D$ and suppose $(\phi-1)D\in \p D$. Then
$\log \phi\in \p\End_{\Z_p}(D)$ is well-defined and it is a
derivation of $D$.
\end{lemma}
\begin{proof}
From the hypothesis  on $\phi$ it follows that $\phi\in 1+\p
\End_{\Z_p}(D)$. Since $D$ is complete in the topology induced by
the filtration  $p^iD$, $\End_{\Z_p}(D)$ is complete in the
topology induced by the filtration  $p^i\End_{\Z_p}(D)$. By the
previous theorem there exists $\log \phi$.

Let
\begin{equation*}
f_n(t)=\sum_{i=1}^n \frac{(-1)^{i+1}(t-1)^i}{i}=\sum_i \alpha_it^i.
\end{equation*}
 Define $d_n=f_n(\phi)$. Note that $\log \phi=\lim_{n\to \infty}d_n$. Then $$d_n(ab)=\sum_i\alpha_i\phi^i(a)\phi^i(b)=\sum_{j,k}\beta_{j,k}(\phi-1)^j(a)(\phi-1)^k(b),$$
 where $\beta_{j,k}$ are coefficients obtained from the following equality:
$$f_n(ts)=\sum_{j,k}\beta_{j,k}(t-1)^j(s-1)^k.$$
From the definition of $f_n$, it follows that
$$f_n(ts)-f_n(t)-f_n(s)=\sum_{j+k>n}\beta_{j,k}(t-1)^j(s-1)^k.$$
Hence
$$d_n(ab)-d_n(a)b-ad_n(b)= \sum_{j+k>n}\beta_{j,k}(\phi-1)^j(a)(\phi-1)^k(b)$$
and so
$$(\log \phi)(ab)-(\log\phi)(a)b-a(\log\phi)(b)=\lim_{n\to \infty}d_n(ab)-ad_n(b)-d_n(a)b=0.$$
Hence $\log \phi$ is a derivation.
\end{proof}
Using the similar argument we can prove the next lemma
\begin{lemma}\label{aut}
Let $D$ be an associative $\Z_p$-algebra without $\Z_p$-torsion and suppose that $\cap_i p^iD=0$. Assume that $D$ is complete in the topology induced by the filtration  $p^iD$. Let $\phi$ be a $\Z_p$-derivation of $D$ and suppose $\phi(D)\in \p D$. Then $\exp(\phi)\in 1+\p\End_{\Z_p}(D)$ is well-defined and it is an automorphism of $D$.
\end{lemma}

Let $R$ be a commutative Noetherian local pro-$p$ domain of
characteristic 0. We  use the notation of the previous section.
 Suppose that $\x\in G(\p R)$.  Then it is clear that the automorphism $\rho_{\x}$
  satisfies the condition: $(\rho_{\x}-1)A\in \p A$,
  whence $\rho_{\x}\in 1+\p\End_R(A)$. By Lemma \ref{derivation},
  we have the well-defined derivation $\log (\rho_{\x})\in \p\End_R(A)$.
  Since $\lambda_{\y}$ and $\rho_{\x}$ commute, $\log (\rho_{\x})\in \p\LL=\p\Der_l(A)$.

From the equality $\exp (\log (\rho_x))=\rho_x$, we obtain,  using
Theorem \ref{explog}(iii), that $G(\p R)$ can be embedded into
$\Gamma(\p\LL)$. In fact, we can prove more:
\begin{theorem}\label{isom}
The groups $\{\rho_{\x}|\x \in G(\p R)\}$ and group $\{e^a|a\in \p \LL\}$ coincide. In particular, $G(\p R)\cong \Gamma(\p\LL)$.
\end{theorem}
\begin{proof}
Let $V$ be the group of $R$-automorphisms of $A$, commuting with all $\lambda_{\x}$, $\x\in G$. We will show that $V=\{\rho_{\x}|\x \in G\}$.

Let $v\in V$ and put $\x=(v(x_1)(\e),\ldots,v(x_n)(\e))$. It is easy to see that if $f\in A$, then $v(f)(\e)=f(\x)=\rho_{\x}f(e)$. Hence if $\y\in G$, we have
$$v(f)(\y)=\lambda_{\y^{-1}}\circ v(f)(\e)=v\circ\lambda_{\y^{-1}}(f)(\e)=\rho_{\x}\circ\lambda_{\y^{-1}}(f)(\e)=\rho_{\x}(f)(\y).$$
It implies that $v=\rho_{\x}$.

Now, if $w\in \p\LL$, then by Lemma \ref{aut}, $e^w$ is an automorphism of $A$. It is clear that $e^w\in V$. Since $e^w(x_i)\in x_i+\p A$, we obtain from the previous paragraph that $e^w=\rho_{\x}$ for some $\x \in G(\p R)$. This finishes the proof.
\end{proof}
\begin{corollary}\label{linear}
If the Lie $R$-algebra $\LL=\LL(G)$ can be embedded in
$\End_R(W)^{(-)}$ for some finitely generated $R$-torsion-free
$R$-module $W$, then $G(\p R)$ can be embedded as a closed
subgroup in $\Aut_R(W)$.
\end{corollary}
\begin{proof}
Suppose $\LL$ is a $R$-subalgebra of $\End_R(W)^{(-)}$. By Theorem
\ref{explog}, $H=\{e^a|a\in\p\LL\}$ is a group  isomorphic to
$\Gamma(\p\LL)$ and, whence, by the previous theorem to $G(\p R)$.
Note that $H$ is  a closed subgroup of $\Aut_R(W)$ because the
exponential map is continuous on $\p\End_R(W)$ and $\p \LL$ is
compact.
\end{proof}

\section{Soluble radical of a finitely generated $R$-standard
group}\label{nilprad} Let $R$ be a noetherian commutative domain,
$D$ its field of fractions and $\LL$ an $R$-Lie algebra which is a
finitely generated free $R$-module. We call $\LL$ for short an
{\bf $R$-lattice}. Put $\LL_D=D\otimes_R \LL $. $\LL_D$ is a
finite dimensional $D$-Lie algebra. In the following $R_s(\LL_D)$
will denote the soluble radical of $\LL_D$ and $R_n(\LL_D)$ will
denote the nilpotent radical of $\LL_D$.
The purpose of this section is
the next result:
\begin{theorem}\label{finitegen}
Let $R$ be a commutative Noetherian local pro-$p$ domain of
characteristic 0 and Krull dimension greater than 1  and $D$ its
field of quotients. Let $G$ be a finitely generated $R$-standard
group and $\LL=\LL(G)$ its Lie algebra. Then $R_s(\LL_D)$ is
nilpotent.
\end{theorem}
We use the notation of the previous section. From Theorem
\ref{isom} we know that if $\x\in G(\p R)$, then $\rho_\x=e^a$ for
some $a\in \p\LL(G)$. Recall that $I$ is an ideal of $A$ generated
by $x_1,\cdots,x_d$. The conjugation by $\x$ which send $f(\y)\in
I$ to $f(\x^{-1}\y\x)\in I$ is  the map $\lambda_\x\circ \rho_\x$.
We need an auxiliary lemma.
\begin{lemma}Let $G$ be a finitely generated
$R$-standard group and $\x\in G(\p R)$. Let $a\in \p\LL(G) $ be
such that $\rho_\x=e^a$. Suppose $\lambda_\x\circ \rho_\x$ acts as
an unipotent automorphism on $(I/I^2)$. Then $a\in R_n(\LL(G)_D)$.
\end{lemma}
\begin{proof}
Let $f\in I$ and $w\in \LL=\LL(G)$. We have
$$(\rho_\x^{-1}\circ w\circ
\rho_\x)(f)(\e)=(w\circ
\rho_\x)(f)(\x^{-1})=(w\circ\lambda_\x\circ \rho_\x)(f)(e).$$
 Note
that if a linear automorphism acts unipotently on $V$, it acts
also unipotently on $V^*$. Since $\LL$ can be identified with
$(I/I^2)^*$ (see Remark \ref{ident}), we have that the
automorphism of $\LL_D$ defined as $\tau(w)= \rho_\x^{-1}\circ
w\circ \rho_\x$ is unipotent.

In order to prove that $\ad a$ is nilpotent, we should to show
that every $\ad a$-invariant subspace $W$ of $\LL_D$ has a nonzero
element $w$ such that $\ad a(w)=0$. So, let $W$ be an $\ad
a$-invariant subspace of $\LL_D$. Since we have
$$\tau(w)=\rho_\x^{-1}\circ w\circ \rho_\x=e^{-a}\circ w\circ e^{a}=\sum_{i=0}^{\infty}\frac{(\ad a)^i(w)}{i!}=e^{\ad a}(w),$$
 $W$ is also $\tau$-invariant. Hence there exists $0\ne w\in W$
 such that $\tau(w)=w$. Then $\ad a(w)=0$. We conclude that $\ad
 a$ is nilpotent.
\end{proof}
\begin{proof}[Proof of Theorem \ref{finitegen}]
We suppose the contrary.  If $v \in R_s(\LL_D)\setminus
R_n(\LL_D)$, then $Dv\cap R_n(\LL_D)=\{ 0 \}$. Put $T=\p\LL\cap
R_s(\LL_D)$ and $T_1=\p\LL\cap Dv$. As in the previous section we
see the elements from  $\LL(G)$ as left invariant derivations of
$A$, and the elements from  $G$ as left invariant automorphisms of
$A$ (so $G$ coincides with $\{\rho_{\x}|\x \in G\}$).

Put $H=\{e^a|a\in T\}$ and $H_1=\{e^a|a\in T_1\}$. By Theorem
\ref{isom}, these two sets are subsets of $G(\p R)=\{\rho_{\x}|\x
\in G(\p R)\}=\{e^a|a\in \p\LL(G)\}$. Note also that $H$ and $H_1$
are subgroups of $G$.

If $h=e^a\in H$ and $g\in G$, then $h^g=(e^{\Ad g(a)})$. Hence $H$
is a normal  subgroup of $G$. On the other hand since $T$ is
soluble, $H$ is soluble.

In  \cite[Proposition 5.1]{Ja} we proved  that for some $k$ the
kernel of the  action by conjugation on $I/I^k$ is $Z(G)$. If
$\Omega$ is the algebraic closure of $D$, we can consider $G/Z(G)$
as a subgroup of $\GL_m(\Omega)$, where $m$ is the rank of $I/I^k$
as $R$-module. Let $\bar G$ and $\bar H$ be the Zariski closures
of $G/Z(G)$ and $HZ(G)/Z(G)$ respectively in $\GL_m(\Omega)$. Then
$\bar H$ is a normal soluble subgroup of $\bar G$.  By \cite[Lemma
19.5]{Hu}, $[\bar G^\circ,\bar G^\circ]\cap \bar H$ is virtually
unipotent.

 Suppose $h=e^a$, $a\in \p \LL$ acts as an
unipotent automorphism on $I/I^k$. Then by the previous lemma,
$a\in R_n(\LL)$. This implies that $H_1Z(G)/Z(G)$ does not have
non-trivial unipotent elements and so $[\bar G^\circ,\bar
G^\circ]\cap H_1Z(G)/Z(G)$ is finite. Since $Dv\cap R_n(\LL_D)=\{
0 \}$, we have $H_1\cap Z(G)=\{1\}$, whence $[\bar G^\circ,\bar
G^\circ]\cap H_1$ is finite.

On the other hand, since $G$ is finitely generated, $G^\circ=G\cap
\bar G^\circ$ is finitely generated and so $H_1/([ G^\circ,
G^\circ]\cap H_1)$ is abelian of finite rank. We obtain that $H_1$
has finite rank. We have a contradiction because $(R,+)$ can be
embedded in $H_1$ and it is not of finite rank.
\end{proof}

\section{The Weigel theorem}\label{prel}
Let $R$ be a noetherian commutative domain, $K$ its field of
fractions and $\LL$ an $R$-Lie algebra which is a finitely
generated free $R$-module.  $\LL_K$ is a finite dimensional
$K$-Lie algebra. The Ado-Iwasawa theorem states that $\LL_K$ has a
finite dimensional linear representation. The next result shows
that if the soluble radical of $\LL_K$ is nilpotent, then this
representation can satisfy some additional nice properties.
\begin{theorem}(T. Weigel,\cite[Lemma 4.3,Proposition 4.4]{We})
 \label{weigel}Let $R$ be an integrally closed noetherian
commutative domain and $K$ its field of fractions. Assume that
$\LL$ is an $R$-lattice and that soluble radical of $\LL_K$ is
nilpotent. Then there exist a finitely generated $R$-torsion-free
$R$-module $W$ and a faithful $R$-linear representation $\psi
\colon \LL\to \End_R(W)$.
\end{theorem}

\section{Linearity of  groups}\label{final}
In this section we finish the proof of Theorem \ref{main}.
\begin{theorem} Let $R$ be a commutative
Noetherian local pro-$p$ ring of characteristic 0 and $G$ be a
finitely generated  $R$-standard group. Then $G$ is $R$-linear.
\end{theorem}
\begin{proof} The theorem is known in the case when $K\dim R=1$.
So we suppose that $K\dim R> 1$. Let  $\F=\F(\x,\y)$
($\x=(x_1,\ldots,x_n),\y=(y_1,\ldots,y_n)$) be the formal  law
associated with $G$ and  $\phi\colon R\to  S$ a homomorphism from
Theorem \ref{embed} when $a=\p$. Then we can extend this
homomorphism to
$$\phi\colon R[[x_1,\ldots,x_n,y_1,\ldots,y_n,z_1,\ldots, z_n ] ]\to S[[x_1,\ldots,x_n,y_1,\ldots,y_n,z_1,\ldots, z_n ] ]$$
in obvious way. Put $\HH (\x,\y)=\phi(\F(\x,\y))$. We have
$$\HH (\HH (\x,\y),\z)=\phi(\F(\F(\x,\y),\z))=\phi(\F(\x,\F(\y,\z)))= \HH (\x, \HH (\y,\z)).$$
Hence $\HH$ is also a formal group law. Let $H$ be an $S$-standard group associated with $\HH$.

Let  $D$ be the
 ring of quotients of $R$. By Theorem \ref{finitegen}, the radical of $D\otimes_R \LL(G)$
is nilpotent. Let $E$ be the field of quotients of $S$. Then we
have
$$E\otimes _S\LL(H)=E\otimes _S(S\otimes _R\LL(G)),$$
which clearly implies that  the radical of $E\otimes_S \LL(H)$ is
also nilpotent.

 Applying $\phi$,  $G(\m^k)$ is embedded as a
closed subgroup into $H(\p S)$:
$$(x_1,\ldots,x_d)\mapsto (\phi(x_1),\ldots,\phi(x_d)).$$
By Theorem \ref{weigel}, $\LL(H)$ acts faithfully on a finitely
generated $S$-torsion-free module $W$. Hence, by Corollary
\ref{linear}, $G(\m^k)\le H(\p S)$ acts faithfully on $W$. By
Theorem \ref{torsionfreemodule}, $G(\m^k)$ is linear over some
commutative Noetherian local pro-$p$ domain $T$ of characteristic
0 and same Krull dimension as  $R$. By Theorem \ref{tlinear},
$G(\m^k)$ is $R$-linear. Finally, since the index of $G(\m^k)$ in
$G$ is finite, $G$ is also  $R$-linear.
\end{proof}

\end{document}